\documentclass[twoside,11pt]{article}
\usepackage{amsthm}
\usepackage{latexsym}
\usepackage{amssymb}
\usepackage{amsfonts}
\usepackage{amsthm}
\usepackage{amsmath}
\usepackage{amscd}

\usepackage{geometry}
\usepackage{comment}
\usepackage{graphicx}
\usepackage{subcaption}

\begin{document}

\begin{center}
{\LARGE \textbf{A Matricial Aspect of Goldbach's Conjecture}}
\end{center}

\begin{center}
\Large Andrei-Lucian Dr\u{a}goi\footnote [1] {E-mail: dr.dragoi@yahoo.com}
\end{center}

\begin{center}
\Large \textbf{Abstract}
\end{center}

\par This article  is a survey based on our earlier paper (``The ``Vertical'' Generalization of the Binary Goldbach's Conjecture as Applied on ``Iterative'' Primes with (Recursive) Prime Indexes (i-primeths)'' \textbf{\cite{Dragoi}}), a paper in which we have proposed a new generalization of the \underbar{binary/``strong"} Goldbach's Conjecture (\textbf{GC}) briefly called ``the Vertical Goldbach's Conjecture'' (\textbf{VGC}), which is essentially a metaconjecture, as VGC states an infinite number of Goldbach-like conjectures stronger than GC, which all apply on ``iterative'' primes with recursive prime indexes (named ``i-primes''). VGC was discovered by the author of this paper in 2007, after which it was improved and extended (by computational verifications) until the present (2019). VGC distinguishes as a very important ``metaconjecture'' of primes, because it states a new class containing an infinite number of conjectures stronger/stricter than GC. VGC has great potential importance in the optimization of the GC experimental verification (including other possible theoretical and practical applications in mathematics and physics). VGC can be also regarded as a very special self-similar property of the primes distribution. This present survey contains some new results on VGC.



\begin{center}
\textbf {***}
\end{center}

\section*{\quad Introduction}

\quad \quad This paper proposes the generalization of the \underbar{binary (aka ``strong'')} Goldbach's Conjecture (\textbf{GC)} \textbf{\cite{Oliveira1}}, briefly called ``the Vertical (binary) Goldbach's Conjecture'' (\textbf{VGC}), which is essentially a meta-conjecture, as VGC states an infinite number of Goldbach-like conjectures stronger/stricter than GC, which all apply on ``iterative'' primes with recursive prime indexes named ``i-primes'' in this paper. 

\begin{center}
\textbf {***}
\end{center}

\section{Notation and terminology}

\quad \quad Given the simplified notation (used in this article) $p_{x} =p\left(x\right)\ge 2$ (with $x\in \mathbb{N}^*=\left\{1,2,3,...\right\}$) as the x-th prime from the infinite countable primes set $P=\left\{p_{1} \left(=2\right),p_{2} \left(=3\right),p_{3} \left(=5\right),...p_{x} ,...\right\}$,\textit{ the ``i-prime'' concept is the generalization with iteration order }$i\in N=\left\{0,1,2,3,...\right\}$ \textit{of the known ``prime-indexed primes'' (alias ``super-primes'') as a subset of (simple or recursive) primes with (also) prime indexes, with }${}^{i} p_{x} =p^{i+1} (x)=p(p(p...(x)))\ne \left[p(x)\right]^{i+1} $ (the primes indexing bijective function ``p'' applied (i+1)-times on any $x\in \mathbb{N}^*$, which implies a number of just i p-on-p iterations) \textit{being the x-th i-prime, with iteration order }$i\in \mathbb{N}$, \textit{as noted in this paper.} In this notation, simple primes are defined as 0-primes so that ${}^{0} p_{x} =p_{x} =p\left(x\right)$ (and ${}^{0} P=P$): ${}^{1} p_{x} =p_{p_{x} } =p^{2} \left(x\right)=p\left(p\left(x\right)\right)$ are called ``1-primes'' $\left({}^{1} P=\left\{{}^{1} p_{1} ,\, ^{1} p_{2} ,...\right\}\subset P\right)$ (with just one p-on-p iteration), ${}^{2} p_{x} =p_{p_{p_{x} } } =p^{3} \left(x\right)=p\left(p\left(p\left(x\right)\right)\right)$ are called ``2-primes'' $\left({}^{2} P=\left\{{}^{2} p_{1} ,\, ^{2} p_{2} ,...\right\}\subset \, ^{1} P\subset P\right)$ (with just two p-on-p iterations) and so on $\left({}^{i} P\subset P\right)$.

There are a number of (relative recently discovered) Goldbach-like conjectures (\textbf{GLCs}) stronger than GC: \textit{these stronger GLCs (including VGC, defined as a collection of an infinite number of GLCs) are tools that can inspire new strategies of finding a formal proof for GC, as I shall try to argue in this paper}\textbf{. }

VGC distinguishes as a very important metaconjecture of primes (with potential importance in the optimization of the GC experimental verification and other possible useful theoretical and practical applications in mathematics [including cryptography and fractals] and physics [including crystallography and M-Theory]), and a very special self-similar property of the primes distribution (\textbf{PD}). 

GC is specifically reformulated by the author of this article as a special property of PD, in order to emphasize the importance of PD study \textbf {\cite {Diamond, Granville1, Granville2, Liang, Soundararajan}}, which PD has multiple interesting \textit{fractal patterns}.

The \underbar{non-trivial GC} (\textbf{ntGC}) variant (excluding the trivial cases for identical prime pairs $p_{x} =p_{y} $) essentially states that: \textit{PD is sufficiently dense and (sufficiently) uniform, so that} \textit{any natural even number }$2n\, \, (with\, n>3)$\textit{ can be written as the sum of at least one Goldbach partition (}\textbf{\textit{GP}}\textit{) of \underbar{distinct primes} }$p_{x} >p_{y} $, so that $2n=p_{x} +p_{y} $.

\begin{center}
\textbf {***}
\end{center}

\section{Various forms of Goldbach's conjecture and other related conjectures. A Synthesis and A/B Classification of the Main Known Goldbach-like Conjectures}

\quad \quad  \textbf{ntGC redefinition}. Let us define the set of GP matrices $M_{4} =\left(5,3\right)$ (with 2$*$4=8=5+3; a single GP with distinct elements), $M_{5} =\left(7,3\right)$ (with 2$*$5=10=7+3; a single GP with distinct elements), $M_{6} =\left(7,5\right)$ (with 2$*$6=12=7+5; a single GP with distinct elements), $M_{7} =\left(11,3\right)$ (with 2$*$7=14=11+3; a single GP with distinct elements), $M_{8} =\left(\begin{array}{l} {13,3} \\ {11,5} \end{array}\right)$ (with 2$*$8=16=13+3=11+5; only two GPs with distinct elements),{\dots} $M_{n} =\left(\begin{array}{l} {P_{x1} ,P_{y1} \left(<P_{x1} \right)} \\ {P_{x2} ,P_{y2} \left(<P_{x2} \right)} \\ {\, \, ...\, \, ,\, \, ...} \end{array}\right)$, with $\left(\begin{array}{l} {p_{x1} +p_{y1} =2n} \\ {p_{x2} +p_{y2} =2n} \\ {...} \end{array}\right)$, $n>3$ (as imposed by the non-triviality GC condition $p_{x} \ne p_{y} $) plus the GP-non-redundancy condition in $M_{n} $ (which actually implies the non-triviality one) $\left(p_{x} >n\right)\Rightarrow \left(p_{y} <n\right)$$\left(\Rightarrow p_{x} \ne p_{y} \right)$, which is an additional condition imposed to eliminate redundant lines of $M_{n} $ (which may contain the same elements as another ``mirror''-line, but in inversed order):  the  $\left(p_{x} >n\right)$ condition also anticipates the ntGC verification algorithm proposed by VGC, which algorithm starts from the $p_{x(?)} $ (which is the prime closest to, but smaller than $2n$), scans all primes $p_{k} \left(\le p_{x(?)} \right)$ downwards to $p_{2} =3$ and verifies the primality of the differences $d_{k} \left(=2n-p_{k} \right)$. 

Based on $M_{n} $ (with entries in $P=\, ^{0} P$ for any $n>3$) and defining the empty matrix $M_{\emptyset } $(with zero lines, thus also zero columns), ntGC can be restated as:\textbf{ }\textit{for any positive integer}\textbf{\textit{ }}\textit{$n>3$}\textbf{\textit{, $M_{n} \ne M_{\emptyset } $}} (\textbf{ntGC})

\begin{center}
\textbf {*}
\end{center}

\quad \textbf{The ``Goldbach-like Conjecture (GLC)'' definition}. A GLC is defined in this paper as \textit{the combination of GC plus any additional conjectured property of $M_{n} $ (other that $M_{n} \ne M_{\emptyset } $)  for any positive integer $n>L$ }(with $L\ge 3$ being also a non-zero positive integer limit)

\quad \textbf{GLCs classification}. GLCs may be classified in \underbar{two major classes} using a double criterion such as:

\begin{enumerate}
\item  \textbf{Type A GLCs} \textbf{(A-GLCs)} are \textit{those GLCs that claim: (}\textbf{\textit{1}}\textit{) Not only that all }$M_{n} \ne M_{\emptyset } $\textit{  for any }$n>L$\textit{ but also (}\textbf{\textit{2}}\textit{) any other non-trivial accessory property/properties common for all }$M_{n} \left(\ne M_{\emptyset } \right)$\textit{. A specific A-GLC is considered ``authentic'' if the other non-trivial accessory property/properties common for all }$M_{n} \left(\ne M_{\emptyset } \right)$\textit{ (claimed by that A-GLC) isn't/aren't a consequence of the 1${}^{st}$ claim (of the same A-GLC). Authentic (at least conjectured as such) A-GLCs are (have the potential to be) ``stronger''/stricter than ntGC as they essentially claim ``more'' than ntGC does.}

\item  \textbf{Type B GLCs} \textbf{(B-GLCs)} \textit{are those GLCs that claim: no matter if all }$M_{n} \ne M_{\emptyset } $\textit{ or just some }$M_{n} \ne M_{\emptyset } $ \textit{for }$n>L$\textit{, all those }$M_{n} $\textit{ that are yet non-}$M_{\emptyset } $\textit{ (for }$n>L$\textit{) have (an)other non-trivial accessory property/properties common for all $M_{n} \ne M_{\emptyset } $\textit{ (for }$n>L)$. A specific B-GLC is considered authentic if the other non-trivial accessory property/properties common for all }$M_{n} \left(\ne M_{\emptyset } \right)$\textit{ (claimed by that B-GLC for }$n>L$\textit{) isn't/aren't a consequence of the fact that some }$M_{n} \ne M_{\emptyset } $ \textit{for }$n>L$\textit{. Authentic (at least conjectured as such) B-GLCs are ``neutral'' to ntGC (uncertainly ``stronger'' or ``weaker'' conjectures) as they claim ``more'' but also ``less'' than ntGC does (although they may be globally weaker and easier to be formally proved than ntGC).}
\end{enumerate}

\begin{center}
\textbf {*}
\end{center}

\quad Other variants of the (generic) Goldbach Conjecture (\textbf{GC}) and related conjectures include:

\begin{enumerate}
\item  \textit {``Any odd integer n$>$5 can be written as the sum of 3 (possibly identical) primes}''. This is the (weak) Ternary Goldbach's conjecture (\textbf{TGC}), which was rigorously proved by Harald Helfgott in 2013 \textbf{\cite {Helfgott1, Helfgott2}} (a complex proof that is generally accepted as valid until present), so that TGC is already considered a proved theorem, and no longer just a conjecture.

\item  \textit {``Any integer n$>$17 can be written as the the sum of exactly 3 distinct primes}''. This is cited as ``Conjecture 3.2'' by Pakianathan and Winfree in their article and is a conjecture stronger than TGC, but weaker than ntGC.

\item  \textit{``Any odd integer n$>$5 can be written as the sum of a prime and a doubled prime [which is twice of any prime]}''. This is Lemoine's conjecture (\textbf{LC}) \textbf{\cite {Lemoine}} which is sometimes erroneously attributed to Levy H. who pondered LC in 1963 \textbf{\cite {Levy}}. LC is stronger than TGC, but weaker than ntGC. LC has also an extension formulated by Kiltinen J. and Young P. (alias the "refined Lemoine conjecture" \textbf{\cite {Kiltinen}}), which is stronger than LC, but weaker than ntGC and won't be discussed in this article (as this paper mainly focuses on those GLCs stronger than ntGC). 

\end{enumerate}

There are also a number of (relative recently proposed) GLCs stronger than ntGC (and implicitly stronger than TGC), that can also be synthesized using the $M_{n} $ definition and A/B GLCs classification:  \textit{these stronger GLCs (as VGC also is) are tools that can inspire new strategies of finding a formal proof for ntGC (but also optimizing the algorithms of ntGC empirical verification up to much higher limits than in the present), as we shall try to argue next. }

\begin{enumerate}

\item  \textbf{\underbar{The Goldbach-Knjzek conjecture [GKC]}} \textbf{\cite {Rivera1}} (which is stronger than ntGC) (slightly  reformulated noting with the even number with $2n$ [not with $n$], so that each ):  ``\textit{For any even integer} $2n>4$, \textit{there is at least one prime number }$p$\textit{ [so that] }$\sqrt{2n} <p\le n$\textit{ and }$q=2n-p$\textit{ is also prime [with }$2n=p+q$\textit{ implicitly]}''. GKC can also be reformulated as: ``\textit{every even integer }$2n>4$\textit{ is the sum of at least one pair of primes with at least one prime in the semi-open interval }$\left(\sqrt{2n} ,n\right]$''. GKC can be also formulated using $M_{n} $  such as:

	\begin{enumerate}
	\item  \textbf{Type A formulation variant}: \textit{ ``For any even integer 2n$>$4, $M_{n} \ne M_{\emptyset }$ and $M_{n}$ contains at least one line with one element $p\in \left(\sqrt{2n} ,n\right]$.''}

	\item  \textbf{Type B (neutral) formulation variant}: \textit{``For any even integer 2n$>$4, those $M_{n}$ which are $\ne M_{\emptyset }$ will contain at least one line with one element $p\in \left(\sqrt{2n} ,n\right]$.''}

	\item  \textbf{A non-trivial GKC (ntGKC) version} (which excludes the cases $p=q=n$) is additionally proposed in this paper (with A/B formulations analogous to the GKC variants) and verified up to $2n=10^{10}$: ``\textit{every even integer }$2n>14$\textit{ is the sum of at least one pair of \underbar{distinct primes} with one prime} $p\in \left(\sqrt{2n} ,n\right)$''.
 	\end{enumerate}

\item  \textbf{\underbar{The Goldbach-Knjzek-Rivera conjecture [GKRC]}} \textbf{\cite {Rivera2}}  (which is also obviously stronger than ntGC, but also stronger than GKC for $n\ge 64$) (reformulated): \textit{``}$\forall \, even\, integer\, \, 2n>4,\, $\textit{there is at least one prime number }$p$\textit{ [so that] }$\sqrt{2n} <p<4\sqrt{2n} $ \textit{and }$q=2n-p$ \textit{is also prime [with }$2n=p+q$\textit { implicitly]}''. GKRC can also be reformulated as: ``$\forall \, even\, integer\, \, 2n>4$, $2n$ \textit{is the sum of at least one pair of primes with one element in the double-open interval }$\left(\sqrt{2n} ,4\sqrt{2n} \right)$''. GKRC can be formulated using $M_{n} $ such as:

	\begin{enumerate}
	\item  \textbf {Type A formulation variant}: \textit { ``For any even integer 2n$>$4, $M_{n} \ne M_{\emptyset }$ and $M_{n}$ contains at least one line with one element $p\in \left(\sqrt{2n} ,4\sqrt{2n} \right)$.''}

	\item  \textbf{Type B (neutral) formulation variant}: \textit {``For any even integer 2n$>$4, those $M_{n}$ which are $\ne M_{\emptyset }$ will contain at least one line with one element $p\in \left(\sqrt{2n} ,4\sqrt{2n} \right)$.''}

	\item  \textbf{A non-trivial GKRC (ntGKRC) version} (which excludes the cases $p=q=n$) is additionally proposed in this paper (with A/B formulations analogous to the GKRC variants) and verified up to $2n=10^{10}$: ``\textit{every even integer }$2n>6$\textit{ is the sum of at least one pair of \underbar{distinct primes} with one prime in the open interval }$\left(\sqrt{2n} ,4\sqrt{2n} \right)$''.

	\end{enumerate}

\item  Noting with $g\left(n\right)$ the number of $M_{n} $ lines for each $2n$ in part (in any GLC) (identical to the standard function $g\left(n\right)$ counting the number of non-redundant GPs for each $2n$ tested in ntGC), any other GLC that establishes an additional superior limit of $g\left(n\right)$ (like Woon's GLC \textbf{\cite {Woon}}) can also be considered stronger that ntGC, because ntGC only suggests $g\left(n\right)>0$ for any positive integer $n>3$ (which implies a greater average number of GPs per each $2n$ than the more selective Woon's GLC does).

\end{enumerate}

\section{The main metaconjecture proposed in this paper: The Vertical Goldbach (meta)conjecture (VGC) - The extension and generalization of ntGC as applied on i-primes}

\quad \quad Alternatively noting ${}^{0} p_{x} =p\left(x\right)$ (the x-th 0-prime, equivalent to the x-th prime in the indexed primes set), ${}^{1} p_{x} =p\left(p\left(x\right)\right)$ (the x-th 1-prime), ${}^{2} p_{x} =p\left(p\left(p\left(x\right)\right)\right)$ (the x-th 2-prime), {\dots}${}^{i} p_{x} =p^{i+1} \left(x\right)$ (the x-th i-prime: not to be confused with the exponential $\left[p\left(x\right)\right]^{i+1} \ne p^{i+1} \left(x\right)$), the (main) analytical variant of VGC (\textbf{aVGC}) states that:
\par ``\textit{For any pair of finite positive integers $\left(a,b\right),with\, a\ge b\ge 0$ defining the (recursive) orders of an a-prime $\left({}^{a} p\right)$ and a distinct b-prime }${}^{b} p$\textit{ $\left({}^{a} p\ne \, {}^{b} p,\, \, but\, not\, necessarily\, {}^{a} p>{}^{b} p\, \, when\, a\ne b\right)$ respectively,  there is a finite positive integer limit $L_{\left(a,b\right)} =L_{\left(b,a\right)} \ge 3$ such as, for any (positive) integer $n>L_{\left(a,b\right)} $ there will always exist at least one pair of  finite positive integer indexes $\left(x,y\right)$ so that }${}^{a} p_{x} \ne {}^{b} p_{y} $ \textit{and ${}^{a} p_{x} +{}^{b} p_{y} =2n$.}''
\par More specifically, VGC states/predicts that: \textit{``the 2D matrix $L$ containing all the $L_{\left(a,b\right)} $ limits (organized on columns indexed with ``a'' and lines indexed with ``b'') has a finite positive integer value in any $\left(a,b\right)$ position, without any catastrophic-like infinities, such as (experimentally verified values): $L_{\left(0,0\right)} =3$, $L_{\left(1,0\right)} =3$, $L_{\left(2,0\right)} =2\, \, 564$,  $L_{\left(1,1\right)} =40\, \, 306$, $L_{\left(3,0\right)} =125\, \, 771$, $L_{\left(2,1\right)} =1\, \, 765\, \, 126$, $L_{\left(4,0\right)} =6\, \, 204\, \, 163$, $L_{\left(3,1\right)} =32\, \, 050\, \, 472$, $L_{\left(2,2\right)} =161\, \, 352\, \, 166$, $L_{\left(5,0\right)} =260\, \, 535\, \, 479$, }$L_{\left(4,1\right)} =?(finite)$, $L_{\left(3,2\right)} =?(finite)$, $L_{\left(3,3\right)} =?(finite)$ {\dots}``

\[L=\left[\begin{array}{cccccc} {\left(3\right)} & {\left(3\right)} & {\left(2\, 564\right)} & {\left(125\, \, 771\right)} & {\left(6\, \, 204\, \, 163\right)} & {\left(260\, \, 535\, \, 479\right)} \\ {\left(3\right)} & {\left(40\, \, 306\right)} & {\left(1\, \, 765\, \, 126\right)} & {\left(32\, \, 050\, \, 472\right)} & {\left(?\right)} & {\left(?\right)} \\ {\left(2\, 564\right)} & {\left(1\, \, 765\, \, 126\right)} & {\left(161\, \, 352\, \, 166\right)} & {\left(?\right)} & {\left(?\right)} & {\left(?\right)} \\ {\left(125\, \, 771\right)} & {\left(32\, \, 050\, \, 472\right)} & {\left(?\right)} & {\left(?\right)} & {\left(?\right)} & {\left(?\right)} \\ {\left(6\, \, 204\, \, 163\right)} & {\left(?\right)} & {\left(?\right)} & {\left(?\right)} & {\left(?\right)} & {\left(?\right)} \\ {\left(260\, \, 535\, \, 479\right)} & {\left(?\right)} & {\left(?\right)} & {\left(?\right)} & {\left(?\right)} & {\left(?\right)} \end{array}\right]\] 

\begin{center}
\textbf {*}
\end{center}

\textbf{Important notes on aVGC:}

\begin{enumerate}
\item  VGC is called a ``\underbar{metaconjecture}'' in a specific sense, with the ``metaconjecture'' concept being defined in this paper as \underbar{an finite/infinite set of strongly related (sub)conjectures}. For each $L_{\left(a,b\right)} $ value we have a distinct conjecture/subconjecture which is noted ``VGC(a,b)'' which briefly states that: ``\textit{all even integers $2n>2\cdot L_{\left(a,b\right)} $ can be written as a sum of at least one pair of distinct i-primes ${}^{a} p_{x} \ne {}^{b} p_{y} $ (named a ``vertical'' GP and abbreviated as ``VGP''), with all $L_{\left(a,b\right)} $ (positive integer) limits being the elements of the previously defined $L$ matrix}. 

\item  The initially defined \textbf{$M_{n} $} matrix (recording all unique GPs for each 2n tested with ntGC) can be generalized as a (2D) matrix \textbf{${}^{\left(a,b\right)} M_{n} $} counting all VGPs generated by each VGC(a,b) as tested on each 2n in part.

\item  The \textbf{$g\left(n\right)$} function (which counts the number of GPs generated by ntGC when tested on each 2n in part) can be also generalized as ${}^{\left(a,b\right)} g\left(n\right)$, which counts the number of VGPs from each \textbf{${}^{\left(a,b\right)} M_{n} $}.

\item  It is obvious that each VGC(a,b) will generate a much narrow ``comet'' than another VGC($c$$<$$a$,$d\mathrm{\le}b$). Defining the average ${}^{\left(a,b\right)} g_{av} \left(k\right):\left[n,m\right]=\sum _{k=n}^{m}\frac{{}^{\left(a,b\right)} g\left(k\right)}{\left(k-n\right)+1}  $, the next graph compares the values of ${}^{\left(0,0\right)} g_{av} \left(k\right)$ as ``filtered'' with GKRC (corresponding to ntGKRC) with the values of ${}^{\left(1,0\right)} g_{av} \left(k\right)$ (corresponding to VGC(1.0)) and ${}^{\left(2,0\right)} g_{av} \left(k\right)$ (corresponding to VGC(2.0)) respectively, all averages applied on the closed interval  $\left[3,5000\right]$.

\begin{figure}[h!]
	\centering
  \includegraphics[width=0.65\linewidth]{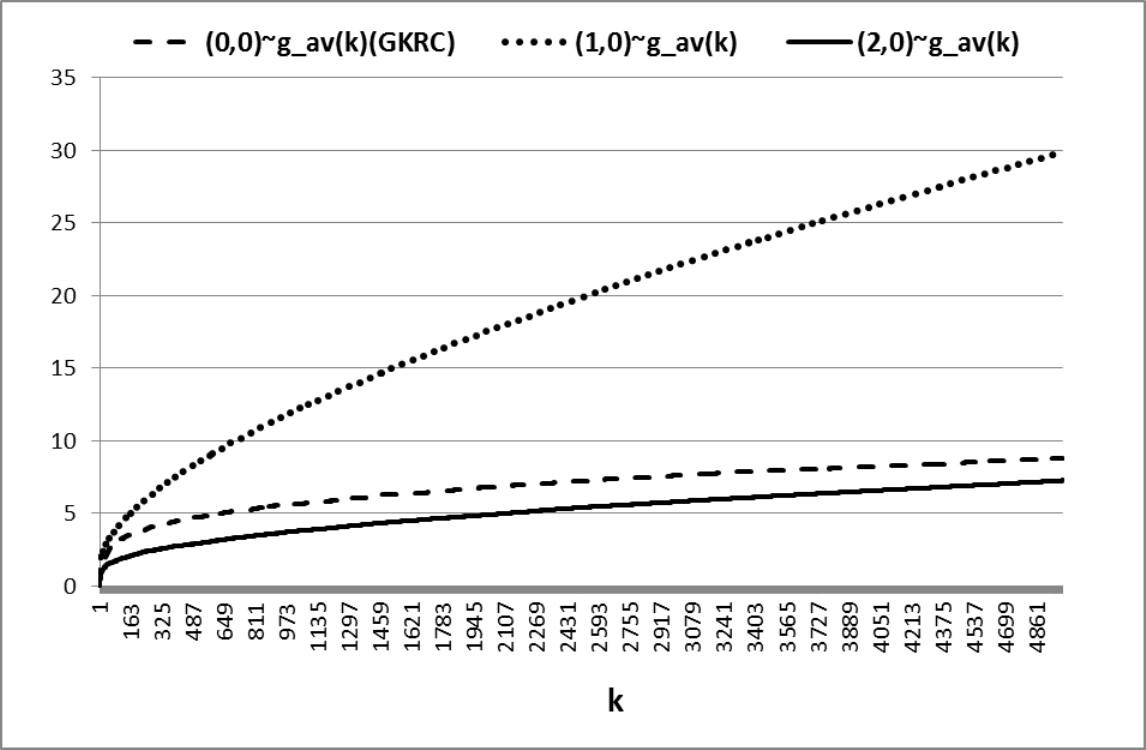}
  \caption{A comparison of the average ${}^{\left(a,b\right)} g_{av} \left(k\right):\left[3,5000\right]$ between VGC(1,0), ntGKRC and VGC(2,0).}
  \label{fig:VGC10_vs_ntGKRC_vs_VGC20}
\end{figure}

\item  aVGC was tested using a Visual C++ software which is shortly named ``\textbf{VGC-SW2}'' \footnote {The source  code of VGC-SW2 (rewritten in Microsoft Studio 2019 - Visual C++ environment as based on an older version called ``VGC-SW1'') that was used to test VGC up to 2n=$10^{10}$ (and also  used for extensive VGC analysis) is available at this URL: drive.google.com/open?id=1egwVzfbc6OyEm9y3-7ISgA3DN0pqTQLS}

\item  The \textbf{$L$} matrix containing only finite values $L_{\left(a,b\right)} $ is essentially a conjectured ``meta-sequence'' of integers and was initially proposed to 
The On-Line Encyclopedia of Integer Sequences (OEIS), but rejected with the main argument that OEIS doesn't accept conjectured meta-sequences and that L matrix was considered ``too ambitious''  \footnote {Review history available at URL: oeis.org/history?seq=A281929\&start=50\ (last page URL); Review history in pdf downloadable format available at URL: dragoii.com/VBGC\_A281929\_OEIS\_rejection\_history.pdf}

\item  The conjectured sequence of all even integers that cannot be expressed as the sum of two distinct 2-prime and 0-prime respectively ${}^{2} p_{x} \ne {}^{0} p_{y} $ was submitted to OEIS, reviewed and approved as \textbf{A282251} sequence (together with VGC cited as reference) \footnote {See URL: oeis.org/A282251; Complete review available at URL: oeis.org/draft/A282251; Review history available at URL: oeis.org/history?seq=A282251}

\item  The conjectured sequence of all even integers that cannot be expressed as the sum of two distinct 1-primes ${}^{1} p_{x} \ne {}^{1} p_{y} $ was also submitted to OEIS, reviewed and approved as \textbf{A316460} sequence (together with VGC cited as reference) \footnote {Sequence entry available at URL: oeis.org/A316460; Complete review available at URL: oeis.org/draft/A316460; Review history available at URL: oeis.org/history?seq=A316460}

\item  Observing the exponential pattern of $L$ elements, a correspondent 2D $G$ matrix is defined as being composed from all $\ln \left(L_{\left(a,b\right)} \right)$ values: 
\[G{\mathop{=}\limits^{\cong }} \left[\begin{array}{cccccc} {\left(1.1\right)} & {\left(1.1\right)} & {\left(7.9\right)} & {\left(11.7\right)} & {\left(15.6\right)} & {\left(19.4\right)} \\ {\left(1.1\right)} & {\left(10.6\right)} & {\left(14.4\right)} & {\left(17.3\right)} & {\left(?\right)} & {\left(?\right)} \\ {\left(7.9\right)} & {\left(14.4\right)} & {\left(18.9\right)} & {\left(?\right)} & {\left(?\right)} & {\left(?\right)} \\ {\left(11.7\right)} & {\left(17.3\right)} & {\left(?\right)} & {\left(?\right)} & {\left(?\right)} & {\left(?\right)} \\ {\left(15.6\right)} & {\left(?\right)} & {\left(?\right)} & {\left(?\right)} & {\left(?\right)} & {\left(?\right)} \\ {\left(19.4\right)} & {\left(?\right)} & {\left(?\right)} & {\left(?\right)} & {\left(?\right)} & {\left(?\right)} \end{array}\right].\] 

\item The $G$ matrix can be graphed as a ``volumetric'' quantized 2D surface in a 3D space: \textbf {see Figure \ref{fig:GMatrixGraph}}

\begin{figure}
	\centering
  \includegraphics[width=0.65\linewidth]{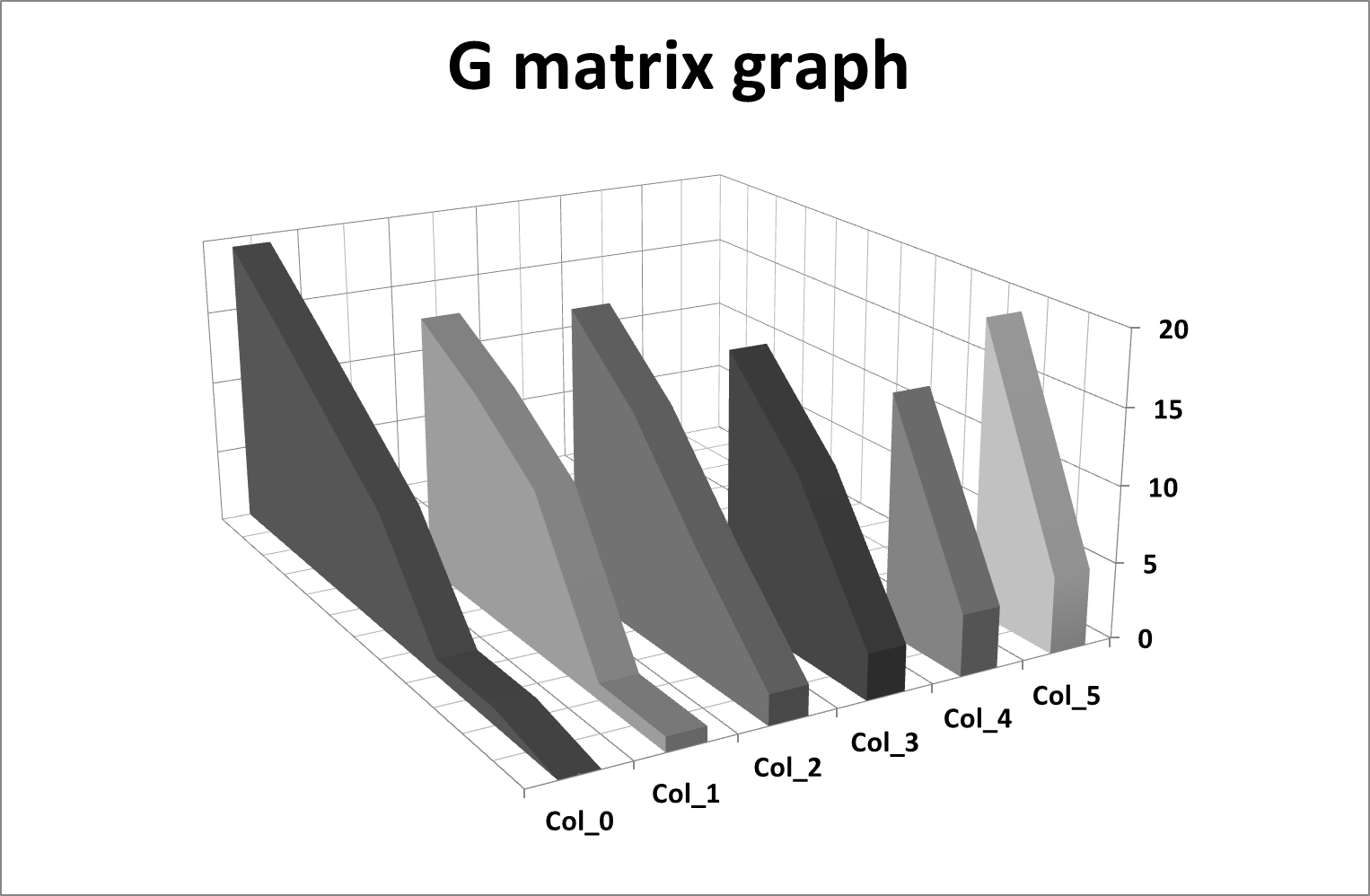}
  \caption{The 3D graph of the $G$ matrix values}
  \label{fig:GMatrixGraph}
\end{figure}

\item  $G$ matrix has a half-``dome''-like graph, apparently with no closed ``depression'' regions, as all elements of $G$ tend to become greater when: moving on the lines from left to right, moving on the columns from up to down, moving on the diagonals, from sides to the center. The elements from each column of $G$ tend to grow almost linearly from up to down (but also on diagonals, from left to center-right and vice versa).

\end{enumerate}

\begin{center}
\textbf {**}
\end{center}

Based on the known (computed) $L_{\left(a,b\right)} $ values (verified up to a superior limit of $2n=10^{10} $), the author also proposes two other (more ambitious) ``inductive'' variants of VGC which aim to predict the magnitude of the unknown $L_{\left(a,b\right)} $ values (marked with ``?'' in L matrix), far beyond the present limit 2n=$10^{10} $: \textbf{see next}.

\begin{center}
\textbf {*}
\end{center}

\textbf{\underbar{The inductive variant of (the meta-conjecture) VGC (iVGC)} proposed in this paper states that}: \textit{``All even positive integers $2n\ge 2\cdot f_{x1} \left(a,b\right)$ and (also) $2n\ge 2\cdot f_{x} {}_{2} \left(a,b\right)$ can be written as the sum of at least one pair of distinct odd primes ${}^{a} p_{x} \ne {}^{b} p_{y} $, with}:
\[f_{x1} \left(a,b\right)=\left\{\begin{array}{l} {2^{(a+1)(b+1)(a+b+2)} \qquad \qquad for\, \left(a=b=0\right)} \\ {2^{[(a+1)(b+1)(a+b+3)/a]-a} \qquad for\, \left(a=b\right)\, \wedge \, \left(a>0\right)} \\ {2^{(a+1)(b+1)(a+b+2)-(a+b-2)} \qquad for\, \left(a\ne b\right)\, \wedge \, \left[\left(a>0\right)\, \vee \, \left(b>0\right)\right]} \end{array}\right. \] 
\textit{and}
\[fx_{2} \left(a,b\right)=\left\{\begin{array}{l} {2^{(a+1)(b+1)(a+b+2)} \qquad \qquad for\, \left(a=b=0\right)} \\ {2^{[(a+1)(b+1)(a+b+3)/a]-2a} \qquad for\, \left(a=b\right)\, \wedge \, \left(a>0\right)} \\ {2^{(a+1)(b+1)(a+b+2)-(a+b-2)} \qquad for\, \left(a\ne b\right)\, \wedge \, \left[\left(a>0\right)\, \vee \, \left(b>0\right)\right]} \end{array}\right. .'' \] 
\textbf{}

\textbf{Important notes on iVGC:}

\begin{enumerate}
\item \textbf{ }$f_{x1} $ has values (organized in a 2D matrix $X_{1} $) that are strictly larger than (but relative closer to) their correspondent$L_{\left(a,b\right)} $ values: 
\[X_{1} {\mathop{=}\limits^{\cong }} \left[\begin{array}{cccccc} {\left(4\right)} & {\left(128\right)} & {\left(4\, \, 096\right)} & {\left(5.2\times 10^{5} \right)} & {\left(2.7\times 10^{8} \right)} & {\left(5.5\times 10^{11} \right)} \\ {\left(128\right)} & {\left(5.2\times 10^{5} \right)} & {\left(5.4\times 10^{8} \right)} & {\left(7\times 10^{13} \right)} & {\left(...\right)} & {\left(...\right)} \\ {\left(4\, \, 096\right)} & {\left(5.4\times 10^{8} \right)} & {\left(7.6\times 10^{8} \right)} & {\left(...\right)} & {\left(...\right)} & {\left(...\right)} \\ {\left(5.2\times 10^{5} \right)} & {\left(7\times 10^{13} \right)} & {\left(...\right)} & {\left(...\right)} & {\left(...\right)} & {\left(...\right)} \\ {\left(2.7\times 10^{8} \right)} & {\left(...\right)} & {\left(...\right)} & {\left(...\right)} & {\left(...\right)} & {\left(...\right)} \\ {\left(5.5\times 10^{11} \right)} & {\left(...\right)} & {\left(...\right)} & {\left(...\right)} & {\left(...\right)} & {\left(...\right)} \end{array}\right]\] 

\begin{enumerate}
\item \textbf{ }$L_{\left(4,1\right)} $ is expected to be smaller than $L_{\left(3,2\right)} $ according to the prediction $f_{x1} \left(4,1\right)<f_{x1} \left(3,2\right)$; obviously $L_{\left(4,1\right)} $ is expected to be larger than $L_{\left(3,1\right)} $ as also according to the prediction $f_{x1} \left(4,1\right)>f_{x1} \left(3,1\right)$; $L_{\left(4,1\right)} $ is ALSO expected to be larger than $L_{\left(3,3\right)} $, according to the prediction $f_{x1} \left(4,1\right)>f_{x1} \left(3,3\right)$; however, $f_{x1} \left(4,1\right)\cong 10^{20} $ surely overestimates $L_{\left(4,1\right)} $, as explained later on in this paper.

\item  \textbf{$f_{x1} \left(3,2\right)\cong 10^{24} $ }surely overestimates $L_{\left(3,2\right)} $, as also explained later on in this paper.

\item  \textbf{$f_{x1} \left(3,3\right)\cong 10^{13} $} almost surely overestimates $L_{\left(3,3\right)} $ over $2n=10^{10} $, as shall also explained later on.

\end{enumerate}

\item  In comparison,$f_{x2} $ has values (organized in a 2D matrix $X_{2} $) that are strictly larger than (but more closer to) their correspondent $L_{\left(a,b\right)} $ values: however, $f_{x2} $ (wrongly) predicts inversed inequalities between some known $L_{\left(a,b\right)} $ values.
\[X_{2} {\mathop{=}\limits^{\cong }} \left[\begin{array}{cccccc} {\left(4\right)} & {\left(128\right)} & {\left(4\, \, 096\right)} & {\left(5.2\times 10^{5} \right)} & {\left(2.7\times 10^{8} \right)} & {\left(5.5\times 10^{11} \right)} \\ {\left(128\right)} & {\left(2.6\times 10^{5} \right)} & {\left(5.4\times 10^{8} \right)} & {\left(7\times 10^{13} \right)} & {\left(...\right)} & {\left(...\right)} \\ {\left(4\, \, 096\right)} & {\left(5.4\times 10^{8} \right)} & {\left(1.8\times 10^{8} \right)} & {\left(...\right)} & {\left(...\right)} & {\left(...\right)} \\ {\left(5.2\times 10^{5} \right)} & {\left(7\times 10^{13} \right)} & {\left(...\right)} & {\left(...\right)} & {\left(...\right)} & {\left(...\right)} \\ {\left(2.7\times 10^{8} \right)} & {\left(...\right)} & {\left(...\right)} & {\left(...\right)} & {\left(...\right)} & {\left(...\right)} \\ {\left(5.5\times 10^{11} \right)} & {\left(...\right)} & {\left(...\right)} & {\left(...\right)} & {\left(...\right)} & {\left(...\right)} \end{array}\right]\] 

\end{enumerate}

\begin{center}
\textbf {**}
\end{center}

The function \textbf{$h\left(a,b\right)=\left\{\begin{array}{l} {2(a+b+1)^{} \qquad for\, \left(a=b=0\right)} \\ {4(a+b)\qquad for\, \left[\left(a=0\right)\wedge \left(b>0\right)\right]\, \vee \, \left[\left(b=0\right)\wedge \left(a>0\right)\right]} \\ {4(a+b+1)^{} \qquad for\, \left(a>0\right)\, \wedge \left(b>0\right)\qquad } \end{array}\right\}$} also generates a 2D matrix $H$ which has values slightly larger than (but very close to) their correspondent $G$ matrix values, such as:
\[H=\left[\begin{array}{cccccc} {\left(2_{>1.1} \right)} & {\left(4_{>1.1} \right)} & {\left(8_{>7.9} \right)} & {\left(12_{>11.7} \right)} & {\left(16_{>15.6} \right)} & {\left(20_{>19.4} \right)} \\ {\left(4_{>1.1} \right)} & {\left(12_{>10.6} \right)} & {\left(16_{>14.4} \right)} & {\left(20_{>17.3} \right)} & {\left(...\right)} & {\left(...\right)} \\ {\left(8_{>7.9} \right)} & {\left(16_{>14.4} \right)} & {\left(20_{>18.9} \right)} & {\left(...\right)} & {\left(...\right)} & {\left(...\right)} \\ {\left(12_{>11.7} \right)} & {\left(20_{>17.3} \right)} & {\left(...\right)} & {\left(...\right)} & {\left(...\right)} & {\left(...\right)} \\ {\left(16_{>15.6} \right)} & {\left(...\right)} & {\left(...\right)} & {\left(...\right)} & {\left(...\right)} & {\left(...\right)} \\ {\left(20_{>19.4} \right)} & {\left(...\right)} & {\left(...\right)} & {\left(...\right)} & {\left(...\right)} & {\left(...\right)} \end{array}\right]\] 

\begin{center}
\end{center}

\par \underbar{A}\textbf{\underbar{ secondary inductive (form of) (the meta-conjecture) VGC (siVGC)} is also proposed in this paper, stating that}: ``\textit{All even positive integers $2m\ge 2\cdot f_{y} (a,b)$ can be written as the sum of at least one pair of distinct odd primes ${}^{a} p_{x} \ne {}^{b} p_{y} $, with $f_{y} (a,b)=e^{h\left(a,b\right)} $}.''
\par
\par \textbf{Important notes on siVGC:}

\begin{enumerate}
\item \textbf{ }$f_{y} \left(6,0\right)\cong 2.65\times 10^{10} $ predicts an $L_{\left(6,0\right)} $ value which is beyond the verification capabilities of VGC-SW2: this hypothesis was also verified with VGC-SW2 and confirmed that $L_{\left(6,0\right)} $ is larger than the limit $2n=10^{10} $. The exception of VGC(6,0) smaller-and-closest to $2n=10^{10} $ is $9\, \, 997\, \, 202\, \, 434=2\times 4\, \, 998\, \, 601\, \, 217$;\textbf{}

\item \textbf{ }$f_{y} \left(7,0\right)\cong 1.45\times 10^{12} $ predicts an $L_{\left(7,0\right)} $ value which is also far beyond the verification capabilities of VGC-SW2.

\item  $h(a,b)$ is based on a simple step-4 ``rule'' applied from the 1-line top-to-down (or from the 1-column left-to-right), so that the unknown elements of $G$ matrix can be alternatively approximated by adding 4 units to the previous element on the same column or line (but not added to the elements from 0-line or 0-column of $G$), such as:
\[G{\mathop{=}\limits^{\cong }} \left[\begin{array}{cccccc} {\left(1.1\right)} & {\left(1.1\right)} & {\left(7.85\right)} & {\left(11.74\right)} & {\left(15.64\right)} & {\left(19.38\right)} \\ {\left(1.1\right)} & {\left(10.6\right)} & {\left(14.38\right)} & {\left(17.28\right)} & {\left({\mathop{\cong }\limits^{?}} 21.28\right)} & {\left(?\right)} \\ {\left(7.85\right)} & {\left(14.38\right)} & {\left(18.9\right)} & {\left({\mathop{\cong }\limits^{?}} 22.9\right)} & {\left(?\right)} & {\left(?\right)} \\ {\left(11.74\right)} & {\left(17.28\right)} & {\left({\mathop{\cong }\limits^{?}} 22.9\right)} & {\left(?\right)} & {\left(?\right)} & {\left(?\right)} \\ {\left(15.64\right)} & {\left({\mathop{\cong }\limits^{?}} 21.28\right)} & {\left(?\right)} & {\left(?\right)} & {\left(?\right)} & {\left(?\right)} \\ {\left(19.38\right)} & {\left(?\right)} & {\left(?\right)} & {\left(?\right)} & {\left(?\right)} & {\left(?\right)} \end{array}\right]\] 

\item  $L_{\left(4,1\right)} =\left(?\right)$ may have a value of $\cong e^{g(3,1)+4} \cong e^{21.28} \cong 1{\rm .7}\times {\rm 10}^{9} $ which value is ALSO under the limit $2n=10^{10} $ and may also be (relatively) verified with VGC-SW2. However, as $L_{\left(4,1\right)} {\mathop{\cong }\limits^{?}} 1{\rm .7}\times {\rm 10}^{9} $ is probably very close to the VGC-SW2 limit $2n=10^{10} $, the subconjecture VGC(4,1) may not be testified by a ``sufficiently'' large gap.

\item \textbf{ }$L_{\left(3,2\right)} {\mathop{\cong }\limits^{?}} e^{g(2,2)+4} \cong e^{22.9} \cong 8{\rm .8}\times {\rm 10}^{9} $ is also predicted (by the step-4 ``rule'') to be under the VGC-SW2 limit $2n=10^{10} $ but ``too close'' to it, so that subconjecture VGC(3,2) may not be testified in VGC-SW2 by a ``sufficiently convincing'' large gap. See the next table.

\end{enumerate}

\begin{table}[h!]
\begin{center}
\caption {The (experimentally) verified values of $L_{\left(a,b\right)}$ (written as exact positive integers) and the estimated maximum $L_{\left(a,b\right)}$ values (using the step-4 ``rule'') smaller than 2n=10$^{10}$ (written in exponential format)}
\label {tab: Table_VGC1}

\begin{tabular}{|p{0.3in}|p{0.7in}|p{0.7in}|p{0.7in}|p{0.7in}|p{0.6in}|p{0.7in}|p{0.6in}|p{0.6in}|} \hline 
\textbf{$L_{\left(a,b\right)} $} & \textbf{0} & \textbf{1} & \textbf{2} & \textbf{3} & \textbf{4} & \textbf{5} & \textbf{6} & \textbf{7            ...} \\ \hline 
\textbf{0} & 3 & 3 & 2,564 & 125,771 & 6,204,163 & 260,535,479 & 1.4E+10 & 7.8E+11 \\ \hline 
\textbf{1} & 3 & 40,306 & 1,765,126 & 32,050,472 & 1.7E+09 & 2.9E+11 & 1.6E+13 & 8.5E+14 \\ \hline 
\textbf{2} & 2,564 & 1,765,126 & 161,352,166 & 8.8E+09 & 2.9E+11 & 1.6E+13 & 8.5E+14 & 4.7E+16 \\ \hline 
\textbf{3} & 125,771 & 32,050,472 & 8.8E+09 & 2.9E+11 & 1.6E+13 & 8.5E+14 & 4.7E+16 & 2.5E+18 \\ \hline 
\textbf{4} & 6,204,163 & 1.7E+09 & 4.8E+11 & 1.6E+13 & 8.5E+14 & 4.7E+16 & 2.5E+18 & 1.4E+20 \\ \hline 
\textbf{5} & 260,535,479 & 2.9E+11 & 2.6E+13 & 8.5E+14 & 4.7E+16 & 2.5E+18 & 1.4E+20 & 7.6E+21 \\ \hline 
\textbf{6} & 1.4E+10 & 1.6E+13 & 1.4E+15 & 4.7E+16 & 2.5E+18 & 1.4E+20 & 7.6E+21 & 4.1E+23 \\ \hline 
\textbf{7\newline ...} & 7.8E+11 & 8.5E+14 & 7.8E+16 & 2.5E+18 & 1.4E+20 & 7.6E+21 & 4.1E+23 & 2.3E+25 \\ \hline 
\end{tabular}
\end{center}
\end{table}

\begin{center}
\textbf {***}
\end{center}

\section{Conclusions on VGC:}

\begin{enumerate}
\item  VGC can be regarded as a ``vertical'' extension and generalization of ntGC as applied on the generalized concept of all subsets of super-primes of any iteration order i, generically named "i-primes" in this paper. VGC has one ``analytical'' variant (aVGC) and two ``inductive'' variants (iVGC and siVGC), which both apply to any subset of primes: 0-primes, 1-primes, 2-primes....

\item  VGC is essentially a \underbar{metaconjecture} in the sense that it actually contains a potential infinite number of subconjectures VGC(a,b), all of them (except VGC(0,0)) stronger/stricter than ntGC.

\begin{enumerate}
\item  \textbf{VGC(0,0)} is equivalent to ntGC.

\item  \textbf{VGC(1,0)} is a GLC stronger and more elegant than ntGC, because it acts on a limit \textbf{$2L_{\left(1,0\right)} =6$} identical to ntGC inferior limit (which is \textbf{$2L_{\left(0,0\right)} =6$}) but is associated with a significantly smaller number of GPs per each even number $2n$.

\item  \textbf{VGC(2,0)} is obviously stronger than VGC(1,0) with much fewer GPs per each $2n$.

\item  \textbf{VGC(1,1}) (anticipated by author's 2007 discovery of VGC(1,0) officially registered in 2012 at OSIM \footnote{For OSIM's product called ``Plicul cu idei'' (``The envelope with ideas'') (by which VGC(1,0) was initially registered by the author) see URL: osim.ro/proprietate-industriala/plicul-cu-idei} is obviously stronger than VGC(1,0) and is equivalent to Bayless-Klyve-Oliveira e Silva Goldbach-like Conjecture (\textbf{BKOS-GLC}) published in Oct. 2013 \textbf {\cite {Bayless}} alias ``Conjecture 9.1'' (rephrased) (tested by these authors up to \textbf{$2n=10^{9} $}): ``\textit{all even integers $2n>\left[2\cdot 40306\left(=2L_{\left(1,1\right)} \right)\right]$ can be expressed as the sum of at least one pair of prime-indexed primes [PIPs] (1-primes }\textbf{\textit{${}^{1} p_{x} $ }}\textit{and ${}^{1} p_{y} $ )}''. This article of Bayless. Klyve and Oliveira (2012, 2013) was based on a previous article by Barnett and Broughan (published in 2009) \textbf{\cite{Broughan}}, but BKOS-GLC was an additional result to this 2009 article. VGC-SW2 was also used to retest and reconfirm VGC(1,1) up to\textbf{$2n=10^{10} $}${}^{\ }$and also helped stating and verifying VGC(a,b) for many other $\left(a,b\right)$ positive integer pairs.

\end{enumerate}

\item  VGC is much ``stronger''/stricter and general than ntGC and proposes a much more rapid and efficient (than the Goldbach-Knjzek-Rivera conjecture [\textbf{GKRC]}) algorithm to find at least one GP for each tested even \textbf{$2n\le 10^{10} $}. 

\begin{enumerate}
\item  VGC is a useful optimized sieve to push forward the limit 4$\cdot$10${}^{18}$ to which GC was verified to hold. 

\item  All VGC(a,b) subconjectures that are distinct from ntGC can be used to produce more rapid algorithms for the experimental verification of ntGC for very large positive integers. A first experiment would be to re-test ntGC up to that limit \textbf{$2n=4\times 10^{18} $} alternatively using various VGC(a,b) and to compare the global times of computing. When verifying ntGC for a very large number $2n$, one can use those aVGC(a,b), iVGC(a,b) or siVGC(a,b) with a minimal positive value for the difference $\left[2n-L_{\left(a,b\right)} \right]$. 

\item  Obviously, VGC(a,b) ``comets'' will tend to narrow progressively with the increase of $a$ and $b$.

\item  VGC is the single known (unified) meta-conjecture of primes and a quite remarkable self-similarity between distinct i-primes and j-primes subsets (with $i>j$), as any i-primes subset (of all primes) is self-similar to the more dense j-primes subset in respect to ntGC, by keeping always finite the limit $2L_{\left(a,b\right)} $ above which all even integers can be written as the sum between two distinct primes ${}^{a} p_{x} \ne {}^{b} p_{y} $. In other words, each of the j-primes subsets behaves as a ``summary of'' any (i$\mathrm{>}$j)-primes set in respect to the ntGC: this is a (quasi)fractal-like GC-related behavior of the primes and primes distribution (\textbf{PD}).

\item  Essentially, VGC conjectures that ntGC is a common property of all i-primes subsets (for any positive integer order i), differing just by the inferior limit $L_{\left(a,b\right)} $ of each VGC(a,b). $L$ matrix is a set of \underbar{critical density thresholds/points} of each i-primes subset in respect to the superset of VGC(a,b) conjectures. All VGC(a,b) comets are self-similar to each other and to the entire VGC superset of comets.

\item  R. G. Batchko has also reported other quasi-fractal structures in the distribution of the prime-indexed primes [\textbf{\cite{Batchko}}: Batchko also used a similar general definition for \underbar{primes with (recursive) prime indexes} (\textbf{PIPs}), alternatively (and briefly) named ``i-primes'' in this paper.

\item  Carlo Cattani and Armando Ciancio also reported a quasi-fractal distribution of primes (including 1-primes) similar to a Cantor set (Cantor dust) by mapping primes and 1-primes into a binary image which visualizes the distribution of 1-primes \textbf{\cite {Cattani}}. 

\item  Obviously, all (i$\mathrm{>}$0)-primes sets $\left({}^{i>0} P\right)$ are subsets of 0-primes set (${}^{0} P=P$) and come in infinite number. There are a potential infinite number of rules/criterions/theorems to extract an infinite number of subsets from $\wp $ (grouped in a family of subsets defined by that specific rule/criterion/theorem: like the Dirichlet's theorem on arithmetic progressions for example). It would be an interesting research subfield of GC to test which are those families (of subsets of primes) that respect ntGC and generate functions with finite values similar to $L_{\left(a,b\right)} $. This potential future research subfield may also help in optimizing the algorithms used in the present for ntGC verification on large numbers. 

\item  It is a quite remarkable fact per se that all ${}^{\left(i>0\right)} P$ subsets have very low densities (when compared to ${}^{0} P=P$), but \underbar{these low densities are sufficiently large} to allow the existence of a matrix $L$ with finite values $L_{\left(a,b\right)}$ for any pair of finites $\left(a,b\right)$. In other (more plastic) words, ntBGC appears to be just a "tree" in the plausibly infinite VGC "wood", which VGC is a spectacular quasi-fractal property of the primes distribution when iteratively applied on itself, generating an infinite number of i-primes subsets that all hold VGC (and ntBGC implicitly) above specific finite integer limits.

\item  A real challenge in the future (concerning VGC) is to calculate the limits $L_{\left(a,b\right)} $ and test/verify other VGC(a,b) subconjectures for large positive integers pairs (a,b), including the pairs $\left(a,b\right)$with relatively large $\left(a-b\right)$ differences.

\end{enumerate}

\end{enumerate}

\begin{center}
\textbf {***}
\end{center}

\section{Potential Applications of VGC:}

\begin{enumerate}
\item \textbf{ }Because the (weak) Ternary Goldbach Conjecture (\textbf{TGC}) is considered a consequence of GC, VGC can be used as a model to also formulate a ``Vertical'' (generalization) of the Ternary Goldbach Conjecture (\textbf{VTGC}) as an analogous consequence of VGC, with a corresponding (potential infinite) meta-sequence of conjectures VTGC(a,b,c), each with an associated limit $L_{\left(a,b,c\right)} $.

\item  VGC can be used to optimize the algorithms of finding very large new primes (i-primes) smaller but closest possible to a chosen (very large) even number $q=2n$: 

\begin{enumerate}
\item  \textbf{Step 1}. One may choose an a-primes subset ${}^{a} P$ and a conjecture VGC(a,b) with positive integer order $b$ chosen so that the known $L_{\left(a,b\right)} $ to be smaller but closest possible to $q=2n$.

\item  \textbf{Step 2}. One may then test only the primality of the differences $d_{x} =q-\left({}^{a} p_{x} \right)$ (starting from $x=1$ to larger positive integer $x$ indexes, in ascending order) which have the potential to be b--primes$\left({}^{b} p_{y} \right)$.

\end{enumerate}

\item  VGC can offer a rule of asymmetric decomposition of Euclidean/non-Euclidean finite/infinite spaces with a finite (positive integer) number of dimensions $d=2n$ into products of pairs of spaces, both with a (positive) i-prime number of dimensions. According to VGC, a finite regular Euclidian/non-Euclidean 2n-space with volume $V_{2n} $ (with n$\mathrm{>}$2) can always be decomposed to a dimensionally asymmetric product of volumes such as: 
\[V_{2n} =V_{\left({}^{a} p_{x} \right)} \times V_{\left({}^{b} p_{y} \right)} =k\cdot \left(r^{\left({}^{a} p_{x} \right)} \times r^{\left({}^{b} p_{y} \right)} \right),with\, \, k=space\, \, volume\, \, specific\, \, constant\] 

\begin{enumerate}
\item  In this way, VGC can also be used in M-Theory to simulate asymmetrical decompositions of 2n-branes (with finite [positive] integer number of dimensions $d=2n$) into products of an ${}^{a} p_{x} $-brane and ${}^{b} p_{y} $ -brane, each brane with a (positive) distinct i-prime number of dimensions.
\end{enumerate}

\item  This type of vertical generalization (generating a meta-conjecture) may be the start of a new research sub-field on primes in which other conjectures may be hypothesized to also have vertical generalizations applied on i-primes. For example, a hypothetical vertical Polignac's conjecture (a ``minus'' version of GC: ``\textit{for any positive even number 2n, there are infinitely many prime gaps of size 2n}'' or ``\textit{there are infinitely many cases of two consecutive prime numbers with even integer difference 2n}'') may speed up the searching algorithms to find very large primes (smaller but closest to a chosen positive even integer superior limit $2n$).

\item  The set of siVGC(a,0) conjectures can also be used to verify much more rapidly (cost/time-efficiently) ntGC, by searching using \textit{only} the subsets ${}^{a} P$ , starting from an ${}^{a} P_{x} $ closest to $2n\ge 2f_{y} (a,0)$ down to ${}^{a} p_{2} $ and testing the primality of $\left(2n-{}^{a} p_{x} \right)$

\end{enumerate}

\begin{center}
\textbf {***}
\end{center}

\section{Addendum: The Short Description of the Software ``VGC-SW2'' Created and Used to Verify VGC}

\quad VGC-SW1 was updated to VGC-SW2 in Microsoft Visual C++ (part of Microsoft Visual Studio 2017 environment).  At first, VGC-SW1 was used to create (and store on hard-disk) a set of ``.bin'' files containing all known i-primes in the double-open interval\textbf{$\left(1,10^{10} \right)$}: \textbf{see the next 
.}

\begin{table}[h!]
\begin{center}
\caption {The files used by the "VGC-SW2" software.}
\label {tab: Table_VGC2}
\begin{tabular}{|p{0.7in}|p{1.6in}|p{1.0in}|p{1.0in}|} 
\hline \textbf{Set of i-primes} & \textbf{File storing the set of i-primes} & \textbf{File dimension on hard-disk (non-archived)} & \textbf{Number of i-primes stored in the file} \\ \hline 
0-primes & 0\_Px\_up\_to\_10 $\mathrm{\hat{}}$ 10.bin & $\mathrm{\sim}$3.55 Gb & {\dots} \\ \hline 
1-primes & 1\_Px\_up\_to\_10 $\mathrm{\hat{}}$ 10.bin & $\mathrm{\sim}$188 Mb & \textbf{24,106,415} \\ \hline 
2-primes & 2\_Px\_up\_to\_10 $\mathrm{\hat{}}$ 10.bin & $\mathrm{\sim}$12 Mb & \textbf{1,513,371} \\ \hline 
3-primes & 3\_Px\_up\_to\_10 $\mathrm{\hat{}}$ 10.bin & $\mathrm{\sim}$900 kb & \textbf{115,127} \\ \hline 
4-primes & 4\_Px\_up\_to\_10 $\mathrm{\hat{}}$ 10.bin & $\mathrm{\sim}$86 kb & \textbf{10,883} \\ \hline 
5-primes & 5\_Px\_up\_to\_10 $\mathrm{\hat{}}$ 10.bin & $\mathrm{\sim}$11 kb & \textbf{1,323} \\ \hline 
6-primes & 6\_Px\_up\_to\_10 $\mathrm{\hat{}}$ 10.bin & $\mathrm{\sim}$2 kb & \textbf{216} \\ \hline 
7-primes & 7\_Px\_up\_to\_10 $\mathrm{\hat{}}$ 10.bin & $\mathrm{\sim}$1 kb & \textbf{47} \\ \hline 
\end{tabular}
\end{center}
\end{table}

\quad For every $\left(a,b\right)$ pair with $a\ge b$, VGC-SW1 verified each ${}^{a} p_{x} \left(>{}^{b} p_{x} \right)$ from the intersection (less dense) set   ${}^{a} P\cap \left(2,2n\ge 6\right)$ (starting from the ${}^{a} p_{x} $ closest to $2n-1$ in descending order): it then verified if  the difference $\left(2n-{}^{a} p_{x} \right)$ is an element in the (more) dense set ${}^{b} P$ by using binary section method. 

\quad VGC-SW1 computed each \textbf{$L_{\left(a,b\right)} $ }limit (with the additional condition \textbf{${}^{a} p_{x} \ne {}^{b} p_{y} $} in at least one GP for any \textbf{$n>L_{\left(a,b\right)} $,} with \textbf{${}^{a} p_{x} +{}^{b} p_{y} =2n$}). 

\noindent The computing time for determining and verifying $L_{\left(2,1\right)} =1\; 765\; 126$ and $L_{\left(2,2\right)} =161\; 352\; 166$ (up to $2n=10^{10} $) was about 30 hours in total. The computing time for determining and verifying $L_{\left(3,0\right)} =125\, \, 771$, \textbf{$L_{\left(4,0\right)} =6\; \; 204\, \, 163$  }and \textbf{$L_{\left(5,0\right)} =260\, \, 535\, \, 479$} was also about 30 hours for each value. The computing time for determining and verifying $L_{\left(3,1\right)} =32\, \, 050\, \, 472$ was a couple of days:\textbf{ }no exceptions found between $2\cdot L_{\left(3,1\right)} $ and $2n=10^{10} $ so that $L_{\left(3,1\right)} $ may be a ``veritable'' last exception of VGC(3,1).

VGC-SW2 was subsequently used for extensive analysis of VGC, including VGC(a,b) comets.

\begin{center}
\textbf {***}
\end{center}

\quad \section* {Acknowledgements}

\quad I would like to express all my sincere gratitude and appreciation to all those who offered me substantial and profound inner motivation for the redaction and completion of this manuscript. I would also like to emphasize my friendship with George Anescu (PhD mathematician and programmer) who helped me verify VGC up to $2n=10^{10} $ by creating VGC-SW1 (version 1), which was a much more rapid (by also using parallel programming) and robust software in Visual C++ for this purpose, as an alternative to my first (relatively slow) software created in Visual Basic (initially used to verify VGC up to $2n=10^{9} $ only). In the meantime, I've also created VGC-SW2 (the 2${}^{nd}$ version of VGC-SW, which is based on the VGC-SW1) for extensive statistical analysis of VGC (including VGC comets and other conjectures related to VGC).

\quad My special thanks to professor Toma Albu who had the patience to read my article and the generosity to help and advise me in correcting and creating an essentialized latex version of this paper. Also my sincere gratitude to professor \c{S}erban-Valentin Str\u{a}til\u{a} who adviced me on the first special case VGC(1,0) discovered in 2007 and urged me to look for a more general conjecture based on VGC(1,0).

\begin{thebibliography}{1}

\bibitem{Batchko} 
Batchko, R. G. (2014). \textit {A prime fractal and global quasi-self-similar structure in the distribution of prime-indexed primes}, arXiv:1405.2900v2 [math.GM], 2014, 36 pp.

\bibitem{Bayless} 
Bayless, J.; Klyve, D. and Oliveira e Silva, T (2012, 2013).  \textit {New bounds and computations on prime-indexed primes},  Integers: Annual. 2013;13:17, 23 pp.

\bibitem{Broughan} 
Broughan, K. A and Ross Barnett, A. (2009). \textit {On the subsequence of primes having prime subscripts}. (10 pages). Article no. 09.2.3 from the Journal of Integer Sequences. 2009;12.

\bibitem{Cai1} 
Cai, Y.C. (2002). \textit {Chen's Theorem with Small Primes}. Acta Mathematica Sinica. 2002;18(3):597--604.

\bibitem{Cai2} 
Cai, Y.C. (2008). \textit {On Chen's theorem (II)}. Journal of Number Theory. 2008;128(5):1336--1357.

\bibitem{Cattani} 
Cattani. C. and Ciancio, A (2016). ``On the fractal distribution of primes and prime-indexed primes by the binary image analysis'', Physica A: Statistical Mechanics and its Applications; 2016. DOI: 10.1016/j.physa.2016.05.013.

\bibitem{Chen1} 
Chen, J.R. (1996). \textit {On the representation of a large even integer as the sum of a prime and the product of at most two primes}. Kexue Tongbao. 1966;11(9):385--386.

\bibitem{Chen2} 
Chen, J.R. (1973). \textit {On the representation of a larger even integer as the sum of a prime and the product of at most two primes}. Scientia Sinica. 1973;16:157--176.

\bibitem{Cheng} 
Cheng-Dong, P; Xia-Xi. D.; Yuan, W. (1975). \textit {On the representation of every large even integer as a sum of a prime and an almost prime}. Scientia Sinica. 1975;XVIII(5):599--610

\bibitem{Diamond} 
Diamond, H.G (1982). \textit {Elementary methods in the study of the distribution of prime numbers}. Bull. Amer. Math. Soc. (N.S.). 1982;7(3):553-589. 

\bibitem{Dragoi} 
Drăgoi, A.L. (2017). \textit {The "Vertical" Generalization of the Binary Goldbach's Conjecture as Applied on "Iterative" Primes with (Recursive) Prime Indexes (i-primeths)}. Journal of Advances in Mathematics and Computer Science [JAMCS] 25(2): 1-32, 2017; JAMCS.36895; ISSN: 2456-9968

\bibitem{Gerstein} 
Gerstein, L.J. (1993). \textit {A reformulation of the Goldbach conjecture}. Mathematics Magazine. 1993;66(1):44-45.

\bibitem{Granville1} 
Granville, A. (September 1993). \textit {Harald Cram\'{e}r and the distribution of prime numbers} (based on a lecture presented on 24th September 1993 at the Cram\'{e}r symposium in Stockholm; 1993.

\bibitem{Granville2} 
Granville, A. (2009). \textit {Different approaches to the distribution of primes}. Milan Journal of Mathematic (2009);78:1--25.

\bibitem{Helfgott1} 
Helfgott, H. A. (2013). \textit {The ternary Goldbach conjecture is true}, arXiv:1312.7748v2 [math.NT], 2014, 79 pp.

\bibitem{Helfgott2} 
Helfgott. H. A (2014). \textit {The ternary Goldbach problem}. Snapshots of modern mathematics from Oberwolfach, No. 3/2014; 2014,2015.

\bibitem{Ikorong} 
Ikorong, G.A.N (2008). \textit {Playing with the twin primes conjecture and the Goldbach conjecture}. Alabama Journal of Mathematics, Spring/Fall. 2008;45-52.

\bibitem{Kiltinen} 
Kiltinen, J; Young P. \textit {Goldbach, Lemoine, and a know/don't know problem}. Mathematics Magazine (Mathematical Association of America). 1984;58(4):195--203.

\bibitem{Lemoine} 
Lemoine, E (1894). \textit  {L'interm\'{e}diare des math\'{e}maticiens}. 1(1894):179: ibid 3. 1896; 151

\bibitem{Levy} 
Levy, H (1963). \textit {On Goldbach's conjecture}. Math. Gaz. 1963; 47(1963):274.

\bibitem{Liang} 
Liang, W; Yan, H.; Zhi-cheng, D (2006). \textit {Fractal in the statistics of Goldbach partition}, arXiv:nlin/0601024v2 [nlin.CD], 2006, 16 pp.

\bibitem{Oliveira1} 
Oliveira e Silva, Tom\'{a}s \& Herzog, Siegfried \& Pardi, Silvio (July 2014). \textit {Empirical verification of the even Goldbach conjecture and computation of prime gaps up to $4\cdot 10^{18}$}. 
Mathematics of Computation. 83. 10.1090/S0025-5718-2013-02787-1.

\bibitem{Polymath} 
Polymath D.H.J. (2014). \textit  {The "bounded gaps between primes" Polymath project - a retrospective},  arXiv:1409.8361v1 [math.HO], 2014, 19 pp.

\bibitem{Rivera1} 
Rivera, C (1999-2001). \textit {Conjecture 22. A stronger version of the Goldbach Conjecture (by Mr. Rudolf Knjzek, from Austria)}', web article from Prime Problems \& Puzzles; 1999-2001.

\bibitem{Rivera2} 
Rivera, C. (1999-2001). \textit {Conjecture 22. A stronger version of the Goldbach Conjecture (by Mr. Rudolf Knjzek, from Austria and narrowed by Rivera C.)}, web article from Prime Problems \& Puzzles; 1999-2001. Not published as an article, but only available online as a cited email at URL: www.primepuzzles.net/conjectures/conj\_022.htm

\bibitem{Soundararajan} 
Soundararajan, K (2006). \textit {The distribution of prime numbers}, arXiv:math/0606408v1 [math.NT], 2006, 20 pp.

\bibitem{Sun1} 
Sun, Z.W. (2013, 2014). Chapter \textit {Problems on combinatorial properties of primes} (19 pages) in \textit {Number Theory: Plowing and Starring Through High Wave Forms: Proceedings of the 7th China-Japan Seminar} (edited by Kaneko M., Kanemitsu S. and Liu J.), 2013;2014:169--188.

\bibitem{Sun2} 
Sun, Z.W. (2014). \textit {Towards the twin prime conjecture}.  A talk given at: NCTS (Hsinchu, Taiwan, August 6, 2014), Northwest University (Xi'an, October 26, 2014) and at Center for Combinatorics, Nankai University (Tianjin, Nov. 3, 2014); 2014.

\bibitem{Woon} 
Woon, M.S.C. (2000). \textit {On partitions of Goldbach's conjecture}, arXiv:math/0010027v2 [math.GM], 2000, 6 pp.

\end{thebibliography}
\end{document}